\RequirePackage[l2tabu,orthodox]{nag}
\pdfoutput=1

\newif\ifarxiv%
\arxivtrue%
\ifarxiv%
  \documentclass[format=acmsmall,timestamp=false,screen,nonacm]{acmart}
  \setcopyright{none}
\else
  \documentclass[format=acmsmall,screen,review]{acmart}
  \setcopyright{none}
  \acmPrice{}
\fi
\acmDOI{}
\acmJournal{TOMS}

\usepackage[utf8]{inputenc}
\usepackage{array}
\usepackage{algorithm}
\usepackage{microtype}
\usepackage{algpseudocode}
\usepackage{amsmath}
\usepackage{booktabs}
\usepackage[english]{babel}
\usepackage{enumitem}
\usepackage{hyperref}
\usepackage{pgfplots}
\usepackage{pgfplotstable}
\usepackage{mathtools}
\usepackage{standalone}
\usepackage{tikz}
\usetikzlibrary{perspective,positioning,calc}
\usepackage{verbatim}
\usepackage{subcaption}
\usepackage{xspace}
\usepackage{csvsimple}
\usepackage[group-separator={,},group-minimum-digits=4]{siunitx}

\usetikzlibrary{patterns}

\pgfplotsset{compat=1.17}

\renewcommand{\div}{\operatorname{div}}

\newcommand{\bfb}{\mathbf{b}}
\newcommand{\bfx}{\mathbf{x}}
\def\Nedelec{N\'ed\'elec\xspace}

\newcommand{\transpose}{\ensuremath{\intercal}}
\newcommand{\hdiv}{\ensuremath{H(\operatorname{div})}\xspace}
\newcommand{\hcurl}{\ensuremath{H(\operatorname{curl})}\xspace}

\usepackage{listings}
\definecolor{DarkBlue}{rgb}{0.00,0.00,0.55}
\definecolor{DarkRed}{rgb}{0.55,0.00,0.00}
\definecolor{DarkGreen}{rgb}{0.00,0.55,0.00}
\definecolor{Bittersweet}{rgb}{1.0, 0.44, 0.37}
\definecolor{Purple}{rgb}{0.5, 0.0, 0.5}

\title{FIAT:\ improving performance and accuracy for high-order finite elements}

\author{Pablo D. Brubeck}
\affiliation{%
  \institution{University of Oxford}
  \department{Mathematical Institute}
  \city{Oxford}
  \country{UK}}
\email{brubeckmarti@maths.ox.ac.uk}

\author{Robert C. Kirby}
\affiliation{%
  \institution{Baylor University}
  \department{Department of Mathematics}
  \streetaddress{1410 S.~4th St.}
  \city{Waco}
  \state{TX}
  \country{USA}}
\email{robert_kirby@baylor.edu}

\author{Fabian Laakmann}
\affiliation{\institution{University of Oxford}
        \department{Mathematical Institute}
	\city{Oxford}
	\country{UK}
}
\email{fabian.laakmann@maths.ox.ac.uk}

\author{Lawrence Mitchell}
\affiliation{\institution{NVIDIA Corporation}
  \city{Santa Clara}
  \country{USA}
}
\email{lmitchell@nvidia.com}

\numberwithin{equation}{section}
\usepackage{cleveref}
\crefname{algorithm}{Algorithm}{Algorithms}
\crefname{figure}{Fig.}{Figs.}
\crefname{table}{Table}{Tables}

\begin{abstract}
FIAT (the FInite element Automatic Tabulator) provides a powerful Python library for the generation and evaluation of finite element basis functions on a reference element.
This release paper describes recent improvements to FIAT aimed at improving its run time and the accuracy and efficiency of code generated using FIAT-provided information.
In the first category, we have greatly streamlined the implementation of orthogonal polynomials out of which finite element bases are built. 
The second category comprises several more advances.
For one, we have built an interface to the \lstinline{recursivenodes} package to enable more accurate Lagrange bases at high order.
We have also implemented integral-type degrees of freedom for $\hdiv$ and $\hcurl$ elements, which match the mathematical definitions of the elements more closely and also avoid loss of accuracy in interpolation.
More fundamentally, we have included families of simplicial quadrature rules that require many fewer quadrature points than the Stroud rules previously used in FIAT.\@
Finally, FIAT now provides support for fast diagonalization methods, which enable fast solution algorithms at very high order.
In each case, we describe the new features in FIAT and illustrate some of the gains obtained through simple numerical tests.
\end{abstract}


\begin{document}
\maketitle

\section{Introduction}
Finite element methods provide a powerful suite of tools for the numerical solution of partial differential equations.
The wide range of mesh-based piecewise polynomial approximating spaces for $H^1$, $\hdiv$, $\hcurl$, and other Sobolev spaces gives finite element methods broad applicability, and the ability to vary the order of approximation enables tradeoffs between computational cost and accuracy obtained.
At the same time, implementing this broad range of approximations at various orders presents a technical challenge to the design of flexible and general software.
FIAT, the FInite element Automatic Tabulator, was first introduced
some two decades ago to provide a general tool for just this purpose \citep{Kirby:2004}.
FIAT is continuously developed on GitHub, and point releases are not provided, but this manuscript highlights development over the past few years, describing the state of FIAT as of mid-2024.
More recent efforts, especially a general facility for macroelements,  will be documented in future manuscripts.

At its core, FIAT works with Ciarlet's abstract definition of a finite element as a triple $(K, P, N)$, where
\begin{itemize}
\item $K \subset \mathbb{R}^d$ is a bounded domain with piecewise smooth boundary.  Typically, $K$ is a simple shape such as a simplex or quadrilateral/hexahedron,
\item $P$ is a finite-dimensional function space defined on the closure of $K$, typically consisting of polynomials or vectors/tensors thereof,
\item $N = {\left\{ n_i \right\}}_{i=1}^{\dim{P}}$ is a basis for the dual space $P^\prime$, called the set of \emph{nodes} or \emph{degrees of freedom}.
\end{itemize}
The \emph{nodal basis} for a finite element is the set ${\left\{ \psi_i \right\}}_{i=1}^{\dim P} \subset P$ such that,
\begin{equation}
n_i(\psi_j) = \delta_{i,j}, \ \ \ 1 \leq i, j \leq \dim P.
\end{equation}
The nodes of a finite element typically consist of functionals such as pointwise evaluation of functions or derivatives at particular points, or certain integral moments of functions on $K$ or its boundary facets, and are chosen to enforce certain kinds of continuity between adjacent elements.

The degrees of freedom of a finite element also define a \emph{nodal interpolant}, which maps functions with sufficient smoothness in some possibly infinite-dimensional space into the finite element space.  This is given by
\begin{equation}
\label{eq:ni}
\mathcal{I}(u) = \sum_{i=1}^{\dim P} n_i(u) \psi_i.
\end{equation}

Many finite elements are defined once on a \emph{reference cell}, $\hat{K}$, such as the unit right triangle, and then mapped to each cell in a finite element mesh.  
So, some computational effort can be put into constructing the nodal basis for a reference element $(\hat{K}, \hat{P}, \hat{N})$, as it can be tabulated once and then mapped repeatedly.
FIAT provides a suite of tools for describing polynomial bases over the reference cell, expressed in terms of orthogonal polynomials, which can be stably evaluated to high-order through recurrence relations, and it also provides a rich set of linear functionals, evaluated numerically.
Given some readily computable basis ${\left\{ \widehat{\phi}_i \right\}}_{i=1}^{\dim P}$ for the polynomial space $P$ (for example, an orthogonal basis), we obtain the reference nodal basis with a generalized Vandermonde-type matrix
\begin{equation}
\label{eq:V}
V_{ij} = \widehat{n}_i \left( \widehat{\phi}_j \right).
\end{equation}
The entries of $V^{-1}$ contain the coefficients expanding the nodal basis in terms of the prime basis.
When $\hat{K} = [0, 1]$, $\{ \widehat{\phi}_j \}$ comprise the monomials, and $\{ \widehat{n}_i \}$ pointwise evaluation, the matrix $V$ is just the classic Vandermonde matrix.

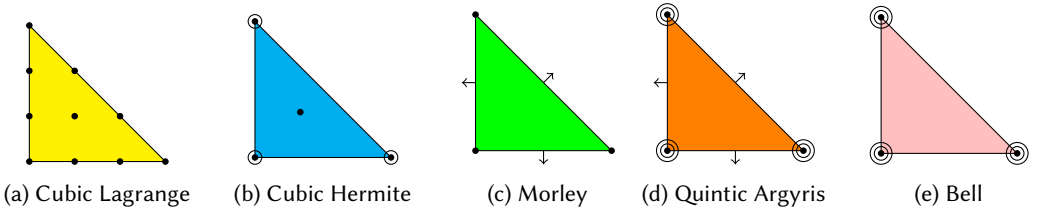
\begin{figure}[htbp]
  \centering
  \begin{subfigure}[t]{0.2\textwidth}
    \centering
    \begin{tikzpicture}[scale=1.8] 
    \draw[fill=yellow] (0,0) -- (1, 0) -- (0, 1) -- cycle;
    \foreach \i in {0, 1, 2, 3} {
      \foreach \j in {0,...,\i}{
      \draw[fill=black] (1-\i/3, \j/3) circle (0.02);
      }
      }
    \end{tikzpicture}
  \caption{Cubic Lagrange\label{lag3}}
  \end{subfigure}
  \hfill
  \begin{subfigure}[t]{0.18\textwidth}
    \centering
  \begin{tikzpicture}[scale=1.8] 
    \draw[fill=cyan] (0,0) -- (1, 0) -- (0, 1) -- cycle;
    \foreach \i/\j in {0/0, 1/0, 0/1}{
      \draw[fill=black] (\i, \j) circle (0.02);
      \draw (\i, \j) circle (0.05);
    }
    \draw[fill=black] (1/3, 1/3) circle (0.02);
  \end{tikzpicture}
  \caption{Cubic Hermite\label{herm3}}
  \end{subfigure}
  \hfill
  \begin{subfigure}[t]{0.18\textwidth} 
    \centering
    \begin{tikzpicture}[scale=1.8]
    \draw[fill=green] (0,0) -- (1, 0) -- (0, 1) -- cycle;
    \foreach \i/\j in {0/0, 1/0, 0/1}{
      \draw[fill=black] (\i, \j) circle (0.02);
    }
    \foreach \i/\j/\n/\t in {0.5/0.0/0.0/-1, 0.5/0.5/0.707/0.707, 0.0/0.5/-1/0}{
      \draw[->] (\i, \j) -- (\i+\n/10, \j+\t/10);
    }
    \end{tikzpicture}
    \caption{Morley\label{morley}}
    \end{subfigure}
  \begin{subfigure}[t]{0.18\textwidth}
    \centering
  \begin{tikzpicture}[scale=1.8] 
    \draw[fill=orange] (0,0) -- (1, 0) -- (0, 1) -- cycle;
    \foreach \i/\j in {0/0, 1/0, 0/1}{
      \draw[fill=black] (\i, \j) circle (0.02);
      \draw (\i, \j) circle (0.05);
      \draw (\i, \j) circle (0.08);
    }
    \foreach \i/\j/\n/\t in {0.5/0.0/0.0/-1, 0.5/0.5/0.707/0.707, 0.0/0.5/-1/0}{
      \draw[->] (\i, \j) -- (\i+\n/10, \j+\t/10);
      }
  \end{tikzpicture}
  \caption{Quintic Argyris\label{arg}}
  \end{subfigure}
  \hfill
  \begin{subfigure}[t]{0.18\textwidth}
    \centering
  \begin{tikzpicture}[scale=1.8] 
    \draw[fill=pink] (0,0) -- (1, 0) -- (0, 1) -- cycle;
    \foreach \i/\j in {0/0, 1/0, 0/1}{
      \draw[fill=black] (\i, \j) circle (0.02);
      \draw (\i, \j) circle (0.05);
      \draw (\i, \j) circle (0.08);
    }
  \end{tikzpicture}
  \caption{Bell\label{bell}}
\end{subfigure}
\caption{Some triangular elements supported by FIAT and their degrees of freedom\label{fig:classicalandzany}}
\end{figure}
Even in its early days, FIAT provided, amongst others, high-order Lagrange elements as well as some of the first realizations of high order Raviart--Thomas \citep{RavTho77a} and \Nedelec{} \citep{nedelec1980mixed} elements (shown later in \cref{fig:hdivels}), giving a major impetus to computing with these ``difficult'' but theoretically powerful elements.
FIAT was an initial component of the FEniCS project \citep{logg2012automated}, and later adopted by the Firedrake project as well \citep{rathgeber2016firedrake, FiredrakeUserManual}.

Over the years, incremental progress has expanded the performance and feature set of FIAT.\@
In \citet{kirby2006optimizing} we showed how to recast much of the element-level computation in terms of matrix multiplication, greatly accelerating the computations.
We also reworked the standard recurrence relations for orthogonal polynomials to avoid coordinate singularities, making it much easier to compute derivatives of basis functions \citep{kirby2010singularity}.
Additionally, FIAT was expanded to include compositional techniques
for defining new elements via tensor product and other operations \citep{mcrae2016automated}.
Building on this, the FInAT project \citep{homolya2017exposing} provides a layer for generating abstract syntax for expressing structured bases that can be exploited in sum-factorization and other optimizations.
This generation of abstract syntax also provided a clean pathway to
incorporate generalized transformations \citep{finat-zany,kirby-zany} for
Argyris \citep{argyris1968tuba} and other elements that do not support
standard reference mappings, such as many of the elements shown in \cref{fig:classicalandzany}.
Later work \citep{aznaran2022transformations,bock2024planar} has included nonstandard transformations for $\hdiv$ elements and nonconforming elements suitable for higher-order equations.
Work on seismic inversion in \citet{spyro} led to the inclusion of
Kong--Mulder--Veldhuizen elements \citep{chin1999higher} and their associated mass lumping quadrature rules.
In addition to the Ciarlet framework, it is also possible to define
the basis functions using a symbolic algebra package such as
sympy \citep{meurer2017sympy}, and we have used this to obtain
serendipity-type elements on quadrilateral and hexahedral elements
\citep{crum2022bringing}.

FIAT was originally a part of the FEniCS project
\citep{logg2012automated}, and is still used in its legacy branch and
in the Firedrake project \citep{FiredrakeUserManual}.
The FEniCSx project has produced a C++ library, basix \citep{scroggs2022basix}, 
that redevelops many features of FIAT within C++, which expedites calling into the basis library at run-time.
Contrasting with this, Firedrake maintains a generative, Python-centric approach.
Recent development effort has greatly improved FIAT's capabilities and performance for high order approximations, and this software release paper documents and demonstrates the resulting gains.

\section{Improved implementation of recurrence relations}\label{sec:rec}
\citet{dubiner1991spectral} introduced an $L^2$-orthogonal polynomial basis over the biunit triangle.
This basis was constructed by mapping the triangle into the biunit square via a singular mapping \citep{duffy1982quadrature}.
There, tensor products of Jacobi polynomials with special weights absorb the Jacobian of the coordinate mapping, and this results in an orthogonal basis on the triangle.
Similar transformations enable the construction of orthogonal bases on tetrahedra, prisms, and pyramids~\citep{karniadakis2005spectral}.

While this tensor product construction allows the application of recurrence relations for univariate polynomials and the construction of sum-factored algorithms, evaluating and especially differentiating the basis at the coordinate singularity requires special care.
Because FIAT frequently requires evaluating and differentiating bases at exactly these points, the work in \citet{kirby2010singularity} provided recurrence relations for the orthogonal polynomials directly on the simplex.
This prevents application of standard Jacobi recurrences but avoids the coordinate singularities altogether.
The key insight for triangles is that the Dubiner polynomials can be written as
\begin{equation}
  D^{p, q}(x, y) = \chi^p(x, y) \psi^{p, q}(y),
\end{equation}
where the factors are in fact polynomials in $x$ and $y$ and are shown to satisfy simple recurrence relations.

Previously, the recurrences from \citet{kirby2010singularity} were evaluated
with a sympy symbol, and each member of the basis was symbolically
differentiated, and then lambdified in order to tabulate the expansion set on a
lattice.  We then obtained a differentiation matrix by solving a Vandermonde
system with the derivative tabulation on the right-hand side.  This operation
was tremendously expensive, and was done every time an element was
instantiated.

In our recent development, we have deemphasized the symbolic approach,
explicitly computing the derivatives of the recurrence relations up to second
order. To achieve this, we unified the recurrences from
\citet{kirby2010singularity} across entity dimensions, so that differentiation
of the three-term recurrence only had to be done in a single place in the code.
Higher order derivatives can be computed through differentiation matrices
mapping coefficients in the orthogonal basis to the coefficients of the partial
derivatives.  We have retained the facility to convert the basis functions to
symbolic form if desired.

\begin{figure}
  \centering
  \begin{subfigure}[t]{0.49\textwidth}
    \centering
  \begin{tikzpicture}[scale=0.8]
    \begin{semilogyaxis}[xlabel={Degree}, ylabel={Time (s)}, xtick=data]
      \addplot[mark=x]  table [x=deg,y=tinit, col sep=comma] {fiat_runtime.Lagrange2.csv};
      \addlegendentry{Init (new)}
      \addplot[mark=x,dashed] table [x=deg,y=teval, col sep=comma] {fiat_runtime.Lagrange2.csv};
      \addlegendentry{Eval (new)}
      \addplot[mark=*] table [x=deg,y=tinit, col sep=comma] {old_fiat_runtime.Lagrange2.csv};
      \addlegendentry{Init (old)}
      \addplot[mark=*, dashed] table [x=deg,y=teval, col sep=comma] {old_fiat_runtime.Lagrange2.csv};
      \addlegendentry{Eval (old)}      
    \end{semilogyaxis}
  \end{tikzpicture}
  \caption{Triangle~\label{lagtri}}
  \end{subfigure}
  \begin{subfigure}[t]{0.49\textwidth}
    \centering
  \begin{tikzpicture}[scale=0.8]
    \begin{semilogyaxis}[xlabel={Degree}, ylabel={Time (s)}, xtick=data]
      \addplot[mark=x]  table [x=deg,y=tinit, col sep=comma] {fiat_runtime.Lagrange3.csv};
      \addlegendentry{Init (new)}
      \addplot[mark=x,dashed] table [x=deg,y=teval, col sep=comma] {fiat_runtime.Lagrange3.csv};
      \addlegendentry{Eval (new)}
      \addplot[mark=*] table [x=deg,y=tinit, col sep=comma] {old_fiat_runtime.Lagrange3.csv};
      \addlegendentry{Init (old)}
      \addplot[mark=*, dashed] table [x=deg,y=teval, col sep=comma] {old_fiat_runtime.Lagrange3.csv};
      \addlegendentry{Eval (old)}      
    \end{semilogyaxis}
  \end{tikzpicture}
  \caption{Tetrahedron~\label{lagtet}}
  \end{subfigure}
  \caption{Time to instantiate Lagrange elements of various orders and tabulate basis functions and first and second derivatives.\label{fiattime}}
\end{figure}
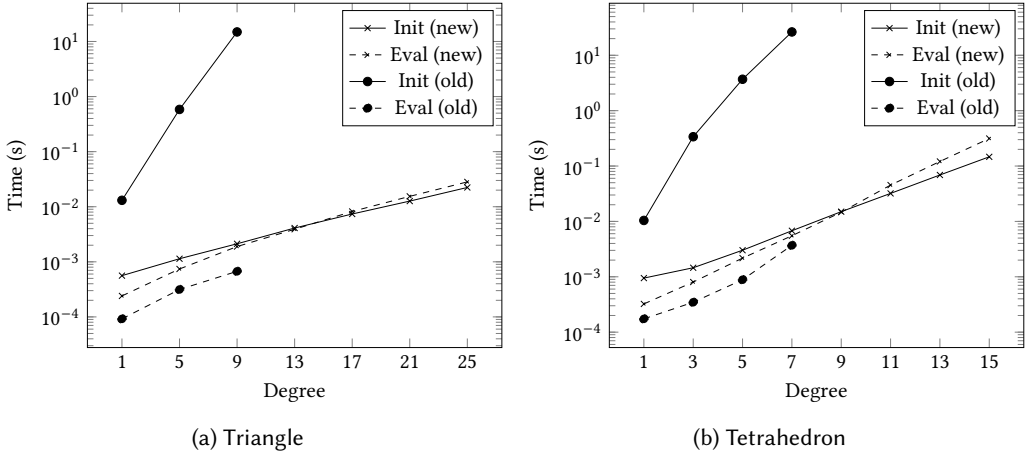

Since FIAT is typically executed as part of a code generation phase, absolute performance is not critical for long-running production simulations since it is executed once during compilation (or during first use of an element at runtime with just-in-time code generation).
Hence, its interface has emphasized flexibility and features over absolute performance.
Early in the modeling process, however, users may frequently re-generate code as they attempt various discretizations and problem formulations and apply them to quite small problems.
Here, the fixed cost of code generation is a more significant part of the user experience.
We document the performance of element instantiation and evaluation before and after our new developments in \cref{fiattime}, where
a few major features appear.
First, the cost of initializing (including forming and inverting the Vandermonde system and sympy conversion) has been massively reduced over the previous version.
Consequently, this allows us to efficiently instantiate the Lagrange basis to much higher order.
Finally, up to second-order derivatives, basis evaluation now evaluates recurrence relations and then applies the inverse Vandermonde matrix rather than applying differentiation matrices.
More time spent in interpreted code leads to a slight uptick in evaluation cost, which is still far outweighed by the gains in construction time.

\section{Better point sets for high-order Lagrange elements}\label{sec:pts}
Historically, FIAT only provided Lagrange bases on the triangle and tetrahedron using points on a uniform point distribution. For triangles, the point set of order $n$ may be described using barycentric coordinates as
\begin{equation}
  \label{eq:bary-point-set}
  \left\{ \left( \frac{i}{n}, \frac{j}{n}, \frac{k}{n} \right): i + j + k = n \right\},
\end{equation}
with an obvious generalization to tetrahedra.
This is sufficient for moderate degrees of approximation, but is well-known to produce difficulties with higher orders of approximation.
For one, equispaced nodes lead to poor conditioning of the
generalized Vandermonde matrix~\cref{eq:V}, even when orthogonal rather than monomial prime bases are used.
This poor conditioning results in roundoff error in the basis tabulation itself, and even if the basis is accurately tabulated, it results in a poor condition number in the overall finite element system.
Moreover, the accuracy of Lagrange interpolation with equispaced nodes is well-known to degrade at high order \citep[Ch.~13]{Trefethen:atap}.

In the univariate case, one
typically uses the Gau\ss--Lobatto points (or Gau\ss--Legendre if
endpoints are not required to be included), which jointly serve as
effective sets of interpolating and quadrature nodes, but the general
coincidence of interpolating points with integration rules is not possible on higher-dimensional simplices.
Various criteria have been proposed for defining suitable node sets on the simplex.  
For example, so-called \emph{Fekete points} maximize a generalized Vandermonde determinant, leading to point distributions in papers such as \citet{chen1995approximate}.
Alternatively, one can start with Fekete points and attempt to optimize the Lebesgue constant of the interpolant \citep{heinrichs2005improved}.
Similarly, \citet{hesthaven1998electrostatics} proposed taking equilibrium positions of electrostatic charges distributed on the simplex.
As processes such as these require a rather expensive optimization phase, they are typically precomputed up to some high order.
Alternatively, it is possible to give explicit constructions (at least in terms of the univariate Lobatto points) of point families that empirically perform quite well.
In this spirit, we point to the warp-blend nodes introduced by \citet{warburton2006explicit}.
More recently, \citet{isaac2020recursive} has given a new family of points
via a recursive, parameter-free construction of simplicial interpolating points that are competitive with the best points in the literature.
He also provides a lightweight Python package called \lstinline{recursivenodes} \citep{recursivenodes} that tabulates these and many other common point families that we utilize within FIAT.\@
Since the warp-blend and recursive constructions rely only on the univariate Lobatto points, they are quite efficient to construct to very high order.

FIAT can now leverage the point distributions available in \lstinline{recursivenodes}.
The \lstinline{Simplex} class, which describes reference geometry and facet connectivity for simplices, has a \lstinline{make_points} method that produces points on the facets.
It now takes an optional keyword argument, \lstinline{variant}, which is then passed into \lstinline{recursivenodes}.
The Lagrange finite element class now also exposes this argument.

We have also plumbed this option into Firedrake, so that users can switch between \lstinline{"equispaced"} and \lstinline{"spectral"} Lagrange variants, the latter being the default.
To demonstrate these new point sets, we performed a small suite of numerical experiments with interpolation and solving the Poisson equation on a coarse mesh with high-degree Lagrange elements. All experiments in this section were performed using Firedrake, see \cref{sec:reproducibility} for code availability.

For interpolation, we divide the biunit square ${[-1,1]}^2$ into two right triangles and interpolate a Runge-type function 
\begin{equation}
  \label{eq:2drunge}
   u(x, y) = \frac{1}{1+(25/2)\left(x^2 + y^2\right)}
\end{equation}
into the $C^0$ Lagrange space with polynomials of degree $k$.
We also divide the biunit cube ${[-1,1]}^3$ into six tetrahedra via the
Freudenthal triangulation \citep{Bey:2000} (shown in
\cref{fig:cube-subdivision}) and again interpolate a Runge-type function
\begin{equation}
  \label{eq:3drunge}
  u(x, y, z) = \frac{1}{1+(25/3)\left(x^2 + y^2 + z^2\right)}.
\end{equation}
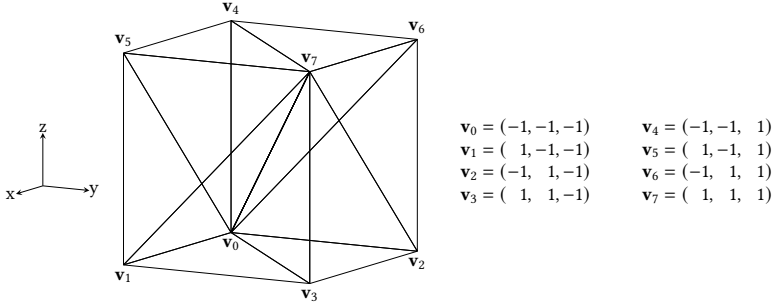
\begin{figure}[htbp]
  \centering
  \begin{tikzpicture}[3d view={120}{10}, scale=2]
  \coordinate [label=below:$\textbf{v}_0$] (v0) at (-1, -1, -1) {};
  \coordinate [label=below:$\textbf{v}_1$] (v1) at (1, -1, -1) {};
  \coordinate [label=below:$\textbf{v}_2$] (v2) at (-1, 1, -1) {};
  \coordinate [label=below:$\textbf{v}_3$] (v3) at (1, 1, -1) {};
  \coordinate [label=above:$\textbf{v}_4$] (v4) at (-1, -1, 1) {};
  \coordinate [label=above:$\textbf{v}_5$] (v5) at (1, -1, 1) {};
  \coordinate [label=above:$\textbf{v}_6$] (v6) at (-1, 1, 1) {};
  \coordinate [label=above:$\textbf{v}_7$] (v7) at (1, 1, 1) {};
  \coordinate (v1v5) at ($(v1)!0.5!(v5)$);
  \coordinate (v3v6) at ($(v3)!0.5!(v6)$);

  \begin{scope}[shift={($(v1v5.west) + (-0.75cm, -0.25cm)$)}]
    \draw[-stealth] (0, 0, 0) -- (0.5, 0, 0) node[pos=1.2] {x};
    \draw[-stealth] (0, 0, 0) -- (0, 0.5, 0) node[pos=1.1] {y};
    \draw[-stealth] (0, 0, 0) -- (0, 0, 0.5) node[pos=1.1] {z};
  \end{scope}
  \draw (v0) -- (v1) (v0) -- (v3) (v0) -- (v7);
  \draw (v1) -- (v3) (v1) -- (v7);
  \draw (v3) -- (v7);

  \draw (v0) -- (v2) (v0) -- (v3) (v0) -- (v7);
  \draw (v2) -- (v3) (v2) -- (v7);
  \draw (v3) -- (v7);

  \draw (v0) -- (v1) (v0) -- (v5) (v0) -- (v7);
  \draw (v1) -- (v5) (v1) -- (v7);
  \draw (v5) -- (v7);

  \draw (v0) -- (v2) (v0) -- (v6) (v0) -- (v7);
  \draw (v2) -- (v6) (v2) -- (v7);
  \draw (v6) -- (v7);

  \draw (v0) -- (v4) (v0) -- (v5) (v0) -- (v7);
  \draw (v4) -- (v5) (v4) -- (v7);
  \draw (v5) -- (v7);

  \draw (v0) -- (v4) (v0) -- (v6) (v0) -- (v7);
  \draw (v4) -- (v6) (v4) -- (v7);
  \draw (v6) -- (v7);

  \node[right=1.5cm of v3v6.east, anchor=west] (A) {
    \begin{tabular}{l}
      $\textbf{v}_0 = (-1, -1, -1)$\\
      $\textbf{v}_1 = (\phantom{-}1, -1, -1)$\\
      $\textbf{v}_2 = (-1, \phantom{-}1, -1)$\\
      $\textbf{v}_3 = (\phantom{-}1, \phantom{-}1, -1)$
    \end{tabular}
  };
  \node[right=1em of A] {
    \begin{tabular}{l}
      $\textbf{v}_4 = (-1, -1, \phantom{-}1)$\\
      $\textbf{v}_5 = (\phantom{-}1, -1, \phantom{-}1)$\\
      $\textbf{v}_6 = (-1, \phantom{-}1, \phantom{-}1)$\\
      $\textbf{v}_7 = (\phantom{-}1, \phantom{-}1, \phantom{-}1)$
    \end{tabular}
  };
\end{tikzpicture}
  \caption{Freudenthal subdivision of the biunit cube.\label{fig:cube-subdivision}}
\end{figure}

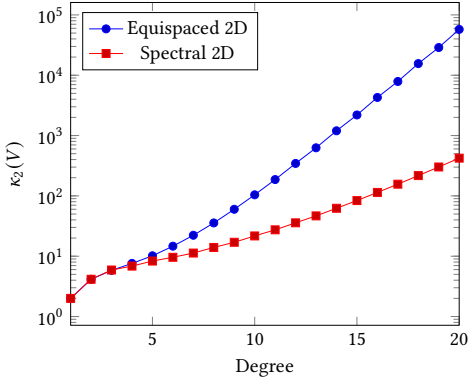
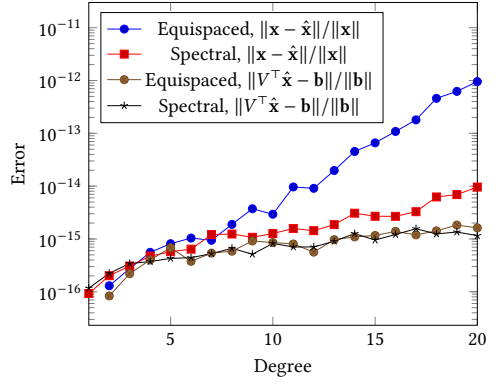
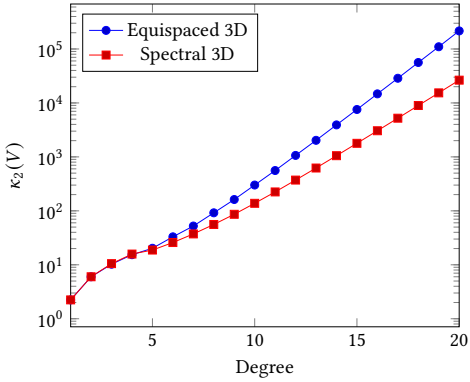
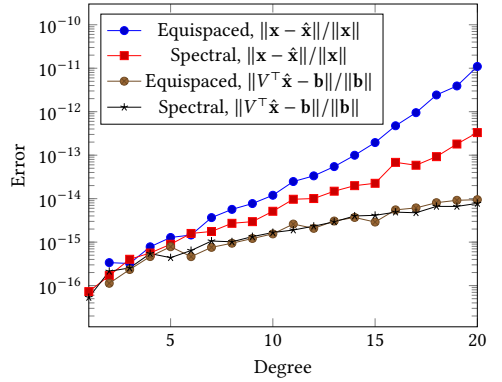
\begin{figure}[htbp]
  \centering
  \begin{subfigure}[c]{0.49\textwidth}
\centering
\begin{tikzpicture}[scale=0.75]    
   \begin{semilogyaxis}[xmin=1,xmax=20,xlabel={Degree}, ylabel= {$\kappa_2(V)$}, legend pos=north west]
        \addplot table [x=deg,y=kappa, col sep=comma] {equispaced.dim2.csv};
        \addlegendentry{Equispaced 2D}
        \addplot table [x=deg,y=kappa, col sep=comma] {gll.dim2.csv};
        \addlegendentry{Spectral 2D} 
   \end{semilogyaxis}
\end{tikzpicture}
\caption{2D conditioning\label{2dcond}}
  \end{subfigure}
\begin{subfigure}[c]{0.49\textwidth}
 \centering
\begin{tikzpicture}[scale=0.75]        
    \begin{semilogyaxis}[xmin=1,xmax=20,xlabel={Degree},ymax=3E-11,
      ylabel={Error}, legend pos=north west]
        \addplot table [x=deg,y=forward, col sep=comma] {equispaced.dim2.csv};
        \addlegendentry{Equispaced, $\| \bfx - \hat{\bfx}\| / \|\bfx\|$}
        \addplot table [x=deg,y=forward, col sep=comma] {gll.dim2.csv};
        \addlegendentry{Spectral, $\| \bfx - \hat{\bfx}\|  / \|\bfx\|$}
        \addplot table [x=deg,y=backward, col sep=comma] {equispaced.dim2.csv};
        \addlegendentry{Equispaced, $\| V^\top \hat{\bfx} - \bfb\| / \|\bfb\|$}
        \addplot table [x=deg,y=backward, col sep=comma] {gll.dim2.csv};
        \addlegendentry{Spectral, $\| V^\top \hat{\bfx} - \bfb\| / \|\bfb\|$}
    \end{semilogyaxis}
\end{tikzpicture}
\caption{2D accuracy\label{2dacc}}
\end{subfigure}\\[\baselineskip]
\begin{subfigure}[c]{0.49\textwidth}
\centering
\begin{tikzpicture}[scale=0.75]    
   \begin{semilogyaxis}[xmin=1,xmax=20,xlabel={Degree}, ylabel= {$\kappa_2(V)$}, legend pos=north west]
        \addplot table [x=deg,y=kappa, col sep=comma] {equispaced.dim3.csv};
        \addlegendentry{Equispaced 3D}
        \addplot table [x=deg,y=kappa, col sep=comma] {gll.dim3.csv};
        \addlegendentry{Spectral 3D} 
   \end{semilogyaxis}
\end{tikzpicture}
\caption{3D conditioning\label{3dcond}}
  \end{subfigure}
\begin{subfigure}[c]{0.49\textwidth}
 \centering
\begin{tikzpicture}[scale=0.75]        
    \begin{semilogyaxis}[xmin=1,xmax=20,xlabel={Degree},ymax=3E-10,
      ylabel={Error}, legend pos=north west]
        \addplot table [x=deg,y=forward, col sep=comma] {equispaced.dim3.csv};
        \addlegendentry{Equispaced, $\| \bfx - \hat{\bfx}\| / \|\bfx\|$}
        \addplot table [x=deg,y=forward, col sep=comma] {gll.dim3.csv};
        \addlegendentry{Spectral, $\| \bfx - \hat{\bfx}\|  / \|\bfx\|$}
        \addplot table [x=deg,y=backward, col sep=comma] {equispaced.dim3.csv};
        \addlegendentry{Equispaced, $\| V^\top \hat{\bfx} - \bfb\| / \|\bfb\|$}
        \addplot table [x=deg,y=backward, col sep=comma] {gll.dim3.csv};
        \addlegendentry{Spectral, $\| V^\top \hat{\bfx} - \bfb\| / \|\bfb\|$}
    \end{semilogyaxis}
\end{tikzpicture}
\caption{3D accuracy\label{3dacc}}
\end{subfigure}
\caption{2-norm conditioning and accuracy of reference Vandermonde inversion with high-degree Lagrange basis using equispaced and recursively-defined GLL nodes.\label{fig:V}}
\end{figure}

FIAT relies heavily on inverting Vandermonde-type matrices based on evaluating orthonormal expansions at the interpolation nodes, so we report on the conditioning and accuracy of these linear systems in \cref{fig:V}.
We see that the recursively-defined GLL-type points in fact give better conditioning than equispaced points, although the gap is somewhat smaller on tetrahedra than triangles.
Also, we report on the accuracy obtained solving systems of the form $V^\top \bfx = \bfb$ by choosing a random solution vector $\bfx$, computing $\bfb$ by matrix-vector multiplication, and then solving the linear system by \lstinline{numpy.linalg.solve}, which interfaces to LAPACK's LU factorization. While the residual norm grows slowly with the degree and does not seem to be affected by the choice of Lagrange points, the error using GLL-type points presents an improvement over the equispaced points, in agreement with the improvement in conditioning.

\begin{figure}[htbp]
  \centering
\begin{subfigure}[c]{0.49\textwidth}
\centering
\begin{tikzpicture}[scale=0.75]
   \begin{semilogyaxis}[xmin=1,xmax=25,xlabel={Degree}, ylabel= {$\|u-\mathcal{I}(u)\|_\infty$}, legend pos=north west]
\addplot table  [x=deg,y=equispaced, col sep=comma] {interp2d.csv};
\addlegendentry{Equispaced 2D}
\addplot table  [x=deg,y=spectral, col sep=comma] {interp2d.csv};
\addlegendentry{Spectral 2D}

\addplot table  [x=deg,y=equispaced, col sep=comma] {interp1d.csv};
\addlegendentry{Equispaced 1D}
\addplot table  [x=deg,y=spectral, col sep=comma] {interp1d.csv};
\addlegendentry{Spectral 1D}
\end{semilogyaxis}
\end{tikzpicture}
\caption{2D\label{interp2d}}
\end{subfigure}
\begin{subfigure}[c]{0.49\textwidth}
\centering
\begin{tikzpicture}[scale=0.75]
   \begin{semilogyaxis}[xmin=1,xmax=15,xlabel={Degree}, ylabel= {$\|u-\mathcal{I}(u)\|_\infty$}, legend pos=north west] 
\addplot table  [x=deg,y=equispaced, col sep=comma] {interp3d.csv};
\addlegendentry{Equispaced 3D}
\addplot table  [x=deg,y=spectral, col sep=comma] {interp3d.csv};
\addlegendentry{Spectral 3D}

\addplot table  [x=deg,y=equispaced, col sep=comma] {interp1d.csv};
\addlegendentry{Equispaced 1D}
\addplot table  [x=deg,y=spectral, col sep=comma] {interp1d.csv};
\addlegendentry{Spectral 1D}
\end{semilogyaxis}
\end{tikzpicture}
\caption{3D\label{interp3d}}
\end{subfigure}
\caption{$L^\infty$ error obtained by interpolating a Runge-type function onto a
coarse mesh with high-degree Lagrange basis using equispaced and
recursively-defined GLL nodes.\label{fig:interpolation-error}}
\end{figure}
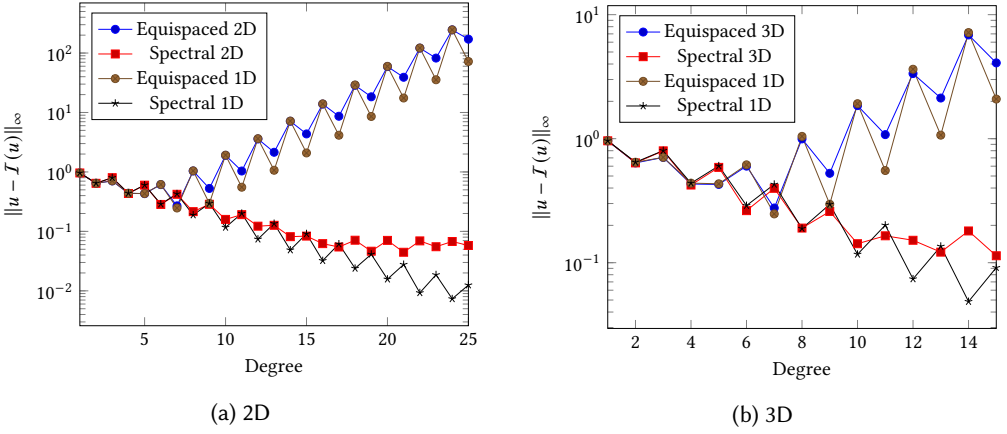
\Cref{fig:interpolation-error} shows the error obtained in interpolating the functions in~\cref{eq:2drunge} and~\cref{eq:3drunge} using both equispaced and the recursively defined GLL-type points.
We observe downward trends (albeit with differences between even and odd degrees) as the degree increases for the recursively-defined nodes.
However, after order 6 or 7, the error using the equispaced points begins to increase dramatically.
In this case, the increasing error results from the exponentially increasing Lebesgue constant (max-norm of the interpolation operator) resulting from large oscillations in the interpolation function.
This generalizes the well-known result of non-convergence of univariate Lagrange interpolation.

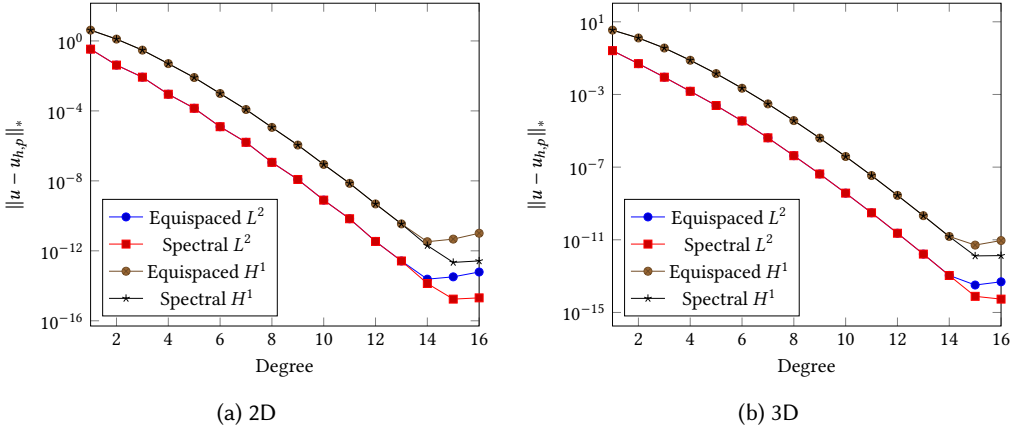
\begin{figure}[htbp]
\begin{subfigure}[c]{0.49\textwidth}
\centering
\begin{tikzpicture}[scale=0.75]
\begin{semilogyaxis}[xmin=1,xmax=16,xlabel={Degree}, ylabel= {$\|u-u_{h,p}\|_*$}, legend pos=south west]
\addplot table  [x=deg,y=equispacedL2, col sep=comma] {poisson2d.csv};
\addlegendentry{Equispaced $L^2$}
\addplot table  [x=deg,y=spectralL2, col sep=comma] {poisson2d.csv};
\addlegendentry{Spectral $L^2$}
\addplot table  [x=deg,y=equispacedH1, col sep=comma] {poisson2d.csv};
\addlegendentry{Equispaced $H^1$}
\addplot table  [x=deg,y=spectralH1, col sep=comma] {poisson2d.csv};
\addlegendentry{Spectral $H^1$}
\end{semilogyaxis}
\end{tikzpicture}
\caption{2D\label{fig:poisson2d}}
\end{subfigure}
\begin{subfigure}[c]{0.49\textwidth}
\centering
\begin{tikzpicture}[scale=0.75]
\begin{semilogyaxis}[xmin=1,xmax=16, xlabel={Degree}, ylabel= {$\|u-u_{h,p}\|_*$}, legend pos=south west]
\addplot table  [x=deg,y=equispacedL2, col sep=comma] {poisson3d.csv};
\addlegendentry{Equispaced $L^2$}
\addplot table  [x=deg,y=spectralL2, col sep=comma] {poisson3d.csv};
\addlegendentry{Spectral $L^2$}
\addplot table  [x=deg,y=equispacedH1, col sep=comma] {poisson3d.csv};
\addlegendentry{Equispaced $H^1$}
\addplot table  [x=deg,y=spectralH1, col sep=comma] {poisson3d.csv};
\addlegendentry{Spectral $H^1$}
\end{semilogyaxis}
\end{tikzpicture}
\caption{3D\label{fig:poisson3d}}
\end{subfigure}
\caption{$L^2$ and $H^1$ errors obtained in solving the Poisson equation on a coarse mesh with high-degree Lagrange basis using equispaced and recursively-defined GLL nodes.\label{fig:poisson-solve-error}}
\end{figure}

Next, we consider solving the Poisson equation using very high orders of
approximation with both equispaced and \lstinline{recursivenodes} GLL
variants, reporting the results in \cref{fig:poisson-solve-error}.
First, we divide the unit square $\Omega = {[0,1]}^2$  into a $4
\times 4$ mesh, each subdivided into two right triangles. We solve
\begin{equation}
  -\Delta u = f,
  \label{eq:poisson}
\end{equation}
subject to homogeneous Dirichlet boundary conditions, with $f$ chosen to ensure that the solution
$u(x, y) = e^{xy} \sin 2\pi x \sin 2 \pi y$.
\Cref{fig:poisson2d} shows that using the recursively defined nodes, the error is decreasing over all polynomial degrees between 1 and 15. In contrast, error starts to increase at degree 13 for the equispaced nodes.
In three dimensions, we similarly divide the unit cube $\Omega = {[0,1]}^3$ into a $4\times 4 \times 4$ grid of cubes, each subdivided into six tetrahedra.
We solve the Poisson equation with homogeneous Dirichlet boundary conditions and chose the forcing function such that the true solution is 
$u(x, y, z) = e^{xyz} \sin 2 \pi x \sin 2 \pi y \sin 2 \pi z$.  
In \cref{fig:poisson3d}, we see error beginning to increase after degree 12 with equispaced points and after degree 13 for the recursively defined nodes, although the increase is less pronounced.

We note that the differences between equispaced and spectral point variants are far less pronounced for a finite element method than for interpolation itself.
Since the Galerkin method finds the \emph{best} energy norm approximation in the approximating space, the methods with two different bases will produce the same mathematical object, up to finite precision effects.
Hence, the differences in \cref{fig:poisson-solve-error} arise from the equispaced points leading to a more ill-conditioned system sensitive to roundoff error.

\section{Integral-type degrees of freedom for \texorpdfstring{{$\hdiv$}}{{H(div)}} and \texorpdfstring{{$\hcurl$}}{{H(curl)}}}\label{sec:intdof}
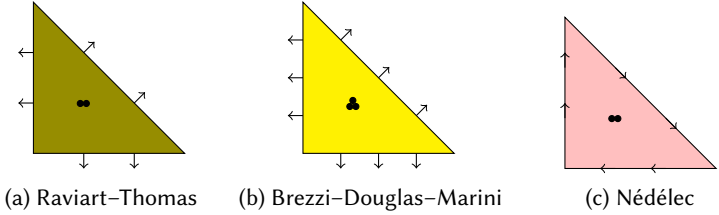
\begin{figure}[htbp]
  \begin{subfigure}[t]{0.25\textwidth}
    \centering
  \begin{tikzpicture}[scale=2.0] 
    \draw[fill=olive] (0,0) -- (1, 0) -- (0, 1) -- cycle;
    \foreach \i/\j/\n/\t in {0.667/0.0/0.0/-1,
                             0.333/0.0/0.0/-1,
                             0.333/0.667/0.707/0.707,
                             0.667/0.333/0.707/0.707,
                             0.0/0.333/-1/0,
                             0.0/0.667/-1/0}{
      \draw[->] (\i, \j) -- (\i+\n/10, \j+\t/10);
    }
   \foreach \i in {0.31,0.35} {
    \draw[fill=black] (\i,0.33) circle (0.02);   
  }
  \end{tikzpicture}
  \caption{Raviart--Thomas\label{rt}}
  \end{subfigure}
  \begin{subfigure}[t]{0.25\textwidth}
    \centering
  \begin{tikzpicture}[scale=2.0] 
    \draw[fill=yellow] (0,0) -- (1, 0) -- (0, 1) -- cycle;
    \foreach \i/\j/\n/\t in {0.75/0.0/0.0/-1,
                             0.5/0.0/0.0/-1,
                             0.25/0.0/0.0/-1,
                             0.25/0.75/0.707/0.707,
                             0.5/0.5/0.707/0.707,
                             0.75/0.25/0.707/0.707,
                             0.0/0.25/-1/0,
                             0.0/0.5/-1/0,
                             0.0/0.75/-1/0}{
      \draw[->] (\i, \j) -- (\i+\n/10, \j+\t/10);
    }
    \foreach \i/\j in {0.31/0.31,0.35/0.31,0.33/0.35} {
    \draw[fill=black] (\i,\j) circle (0.02);
    }
  \end{tikzpicture}
  \caption{Brezzi--Douglas--Marini\label{bdm}}
  \end{subfigure}
  \begin{subfigure}[t]{0.25\textwidth}
    \centering
  \begin{tikzpicture}[scale=2.0] 
    \draw[fill=pink] (0,0) -- (1, 0) -- (0, 1) -- cycle;
    \foreach \i/\j/\n/\t in {0.667/0.0/0.0/-1,
                             0.333/0.0/0.0/-1,
                             0.333/0.667/0.707/0.707,
                             0.667/0.333/0.707/0.707,
                             0.0/0.333/-1/0,
                             0.0/0.667/-1/0}{
      \draw[->] (\i, \j) -- (\i+\t/10, \j-\n/10);
    }
   \foreach \i in {0.31,0.35} {
    \draw[fill=black] (\i,0.33) circle (0.02);   
  }
  \end{tikzpicture}
  \caption{\Nedelec{}\label{ned}}
  \end{subfigure}
  \caption{Some finite elements for $\hdiv$ and
    $\hcurl$.\label{fig:hdivcurl}\label{fig:hdivels}}
\end{figure}

At the time of FIAT's initial development, the excellent theoretical properties $\hdiv$-conforming elements such as Raviart--Thomas \citep{RavTho77a} and Brezzi--Douglas--Marini \citep{brezzi1985two} and the $\hcurl$-conforming elements of \Nedelec{} \citep{nedelec1980mixed,nedelec1986new} (shown in \cref{fig:hdivcurl}), as well as their difficulties in implementation, were well-understood.
Hence, the automatic construction of their bases in FIAT, and our work on form compilation with these elements in \citet{rognes2010efficient} represented a major contribution to finite element computation.
The rapid expansion of finite element exterior calculus \citep{arnold2018finite} has only amplified interest in these elements.

Although these elements are typically defined with integral moments as facet degrees of freedom, $\hdiv$ or $\hcurl$ conformity is also obtained by specifying normal or tangential components at points chosen at the edge.
The initial implementation in FIAT did just this.  
While these degrees of freedom produce suitable finite element approximation, they do not give optimal accuracy for interpolating data into the discrete space \citep[\S 4.2]{laakmann2022augmented}.  
To address this, we have also given implementations of Raviart--Thomas, Brezzi--Douglas--Marini, and first and second kind \Nedelec{} elements using integral moments as facet degrees of freedom.  
Like the selection of optimized point distributions for Lagrange elements, these are also selected through a \lstinline{variant} keyword to the constructor, with choices being \lstinline{"point"} or \lstinline{"integral"}.

Mathematically, the degrees of freedom are defined by integration, but they are implemented by numerical quadrature.
Thus, when computing a nodal interpolant of data by~\cref{eq:ni}, a
quadrature rule that is exact for defining the basis function can incur some numerical error.
It may be necessary in these settings to use a higher order quadrature rule for these degrees of freedom, and the keyword \lstinline{variant="integral(q)"} selects a quadrature rule that has degree \lstinline{q} more than is required for the basis functions.

A higher quadrature degree is  important if one wants to preserve the divergence or curl of a non-polynomial expression during interpolation up to machine precision. If all degrees of freedom are integrated exactly, it holds \citep[Section 2.5]{boffi2013mixed} that
\begin{equation}
	\operatorname{div} \mathcal{I}  (u) =  \Pi \operatorname{div}  u, \quad \text{and} \quad \operatorname{curl} \mathcal{I}   (u) = \Pi \operatorname{curl}  u,
\end{equation}
where $\Pi$ indicates $L^2$ projection into polynomials.  
Enforcing these identities has crucial consequences for structure-preserving finite element approximations, where interpolation is typically applied to source terms or boundary data. In \citet{laakmann2022augmented} it is shown that, both for Maxwell's equations and when solving magnetohydrodynamics, choosing a high enough quadrature degree is necessary to obtain a divergence-free magnetic field on the discrete level;  a property that is important to obtain physically accurate solutions.

\begin{figure}[htbp]
  \centering
  \pgfplotstableread[col sep=comma]{interpRT3d.csv}\datatable%

  \begin{tikzpicture}[scale=0.75]
    \begin{axis}[
      xlabel={Variants},
      ylabel={$\|\operatorname{div} \mathcal{I}(u)\|_{L^2}$},
      xtick=data,
      ybar, log origin=infty,
      xticklabels from table={\datatable}{variants}, 
      xticklabel style={rotate=45, anchor=east},
      legend pos=north east,
      ymode=log, 
      ymin=1e-15, ymax=1, 
      ]
      
      \addplot table [ x expr=\coordindex, y=nref2, ] {\datatable};
      
    \end{axis}
  \end{tikzpicture}
  \caption{$L^2$ norm of $\operatorname{div} \mathcal{I}(u)$ when
    interpolating the non-polynomial divergence-free expression \cref{eq:div-free-u} for different variants and quadrature degrees.\label{fig:interpRT3d}}
\end{figure}
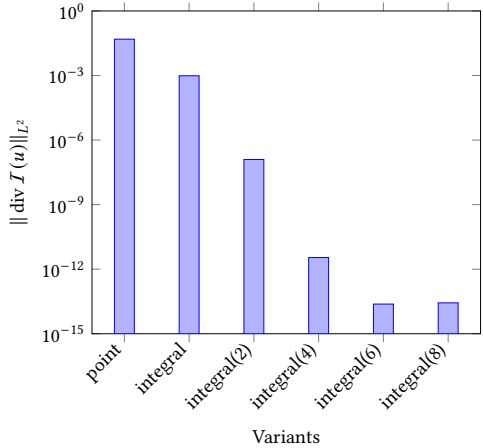

In \cref{fig:interpRT3d}, we show the effect of the quadrature degree
when interpolating the divergence-free function
\begin{equation}
u(x,y,z) = \operatorname{curl}\left({[\sin{(x)} y  \exp{(z)}, \sin{(z)} x y, \cos{(y)}x]}^\transpose\right)
\label{eq:div-free-u}
\end{equation}
into the Raviart--Thomas space of degree two in three dimensions, where $\Omega = {[0,1]}^3$ is divided into an $8\times 8 \times 8$ grid of cubes, each subdivided into six tetrahedra. One can see that increasing the quadrature degree in the numerical integration of the degrees of freedom preserves the divergence more accurately until we reach machine precision for the variant `integral(6)'. While we only report results here for Raviart--Thomas elements in 3D, our experiments show analogous results for Brezzi--Douglas--Marini, and first and second kind \Nedelec{} elements in both two and three dimensions.

We also investigate the convergence order of interpolation for the \lstinline{"point"} and \lstinline{"integral"} variant. For Raviart--Thomas and \Nedelec{} elements of first kind of degree $k$ one expects a convergence order of $k$ in both the $L^2$-norm and $\hdiv$ resp. $\hcurl$ norm. For Brezzi--Douglas--Marini and \Nedelec{} elements of second kind, the convergence order in the $L^2$ norm is increased to $k+1$.

For Raviart--Thomas and \Nedelec{} elements of first kind, we observe a suboptimal convergence order in the $\hdiv$ resp. $\hcurl$ norm for the \lstinline{"point"} variant, while the \lstinline{"integral"} variant produces the expected results. The convergence order for Raviart--Thomas of degree two in three dimensions, where the function $u(x,y,z) = {[\sin{(x)} y  \exp{(z)}, \sin{(z)} x y,
 \cos{(y)}x]}^\transpose$ is interpolated on a base mesh of $2\times 2 \times 2$ cubes over $\Omega = {[0,1]}^3$, each subdivided into six tetrahedra, can be found in \cref{tab:convorderRT3D}. The results for \Nedelec{} elements of first kind are analogous, i.e., only the \lstinline{"integral"} variants shows expected convergence orders. For Brezzi--Douglas--Marini and \Nedelec{} elements of second kind, we observe the expected convergence rates for both variants.
\begin{table}[htbp]
  \centering
    \caption{Convergence order of interpolation into the Raviart--Thomas space of degree two in three dimensions. The \lstinline{"integral"} variant shows the expected convergence order of two in the $\hdiv$ norm, while the \lstinline{"point"} variant shows suboptimal convergence.\label{tab:convorderRT3D}}
    \caption*{variant=\lstinline{"integral"}}
    \begin{tabular}{rllll}
      \toprule
      \# ref & $\|u - \mathcal{I}(u) \|_{L^2}$ & $L^2$ order & $\|u - \mathcal{I}(u) \|_{\hdiv}$ & $\hdiv$ order \\
      \midrule
      0 & 2.99e-02 & -- & 3.50e-02 & -- \\
      1 & 7.54e-03 & 1.99 & 8.85e-03 & 1.98 \\
      2 & 1.89e-03 & 2.00 & 2.22e-03 & 2.00 \\
      3 & 4.73e-04 & 2.00 & 5.55e-04 & 2.00 \\
      \bottomrule
    \end{tabular}
    \par\vspace{\baselineskip}
    \caption*{variant=\lstinline{"point"}}
    \begin{tabular}{rllll}
      \toprule
      \# ref & $\|u - \mathcal{I}(u) \|_{L^2}$ & $L^2$ order & $\|u - \mathcal{I}(u) \|_{\hdiv}$ & $\hdiv$ order \\
      \midrule
      0 & 3.45e-02 & -- & 8.24e-02 & -- \\
      1 & 8.73e-03 & 1.98 & 3.70e-02 & 1.16 \\
      2 & 2.19e-03 & 2.00 & 1.79e-02 & 1.05 \\
      3 & 5.48e-04 & 2.00 & 8.87e-03 & 1.01 \\
      \bottomrule
    \end{tabular}
\end{table}



\section{Expanded suite of efficient quadrature rules}\label{sec:quad}
FIAT originally provided only Stroud conical quadrature \citep{stroud1971approximate} which maps tensor products of Gau\ss--Jacobi rules with appropriate weighting to the unit simplex under the Duffy transform.
While this approach is quite general, the resulting rules are suboptimal at high degree, requiring many more points than the theoretical optimum.
Later, these rules were supplemented by classical hand-coded rules for degrees six and below.

Since use of even moderately high degree finite elements requires quadrature beyond this, we have expanded the suite of optimized quadrature rules.
The numerically-tabulated simplicial quadrature rules of \citet{xiao2010numerical} have become quite popular.
We have included these rules, tabulated up to degree 50 on triangles and degree 15 on tetrahedra.
FIAT's \lstinline{create_quadrature} scheme now selects a low order hand-coded rule where it is better than the Xiao--Gimbutas rule for the same degree, then the Xiao--Gimbutas where available, and Stroud rules as a fall-back for very high order.

\begin{figure}
\begin{subfigure}[c]{0.49\textwidth}
\centering
\begin{tikzpicture}[scale=0.75]
  \begin{semilogyaxis}[xlabel={Degree}, ylabel= {Num points}, legend pos=north west]
\addplot table  [x=deg,y=Stroud, col sep=comma] {quadcard2.csv};
\addlegendentry{Stroud}
\addplot table  [x=deg,y=XG, col sep=comma] {quadcard2.csv};
\addlegendentry{XG}
\end{semilogyaxis}
\end{tikzpicture}
\caption{Triangle\label{quad2}}
\end{subfigure}
\begin{subfigure}[c]{0.49\textwidth}
\centering
\begin{tikzpicture}[scale=0.75]
\begin{semilogyaxis}[xlabel={Degree}, ylabel= {Num points}, legend pos=north west]
\addplot table  [x=deg,y=Stroud, col sep=comma] {quadcard3.csv};
\addlegendentry{Stroud}
\addplot table  [x=deg,y=XG, col sep=comma] {quadcard3.csv};
\addlegendentry{XG}
\end{semilogyaxis}
\end{tikzpicture}
\caption{Tetrahedron\label{quad3}}
\end{subfigure}
\caption{Cardinality of Xiao--Gimbutas versus Stroud conical rules on triangles and tetrahedra.\label{quad}}
\end{figure}
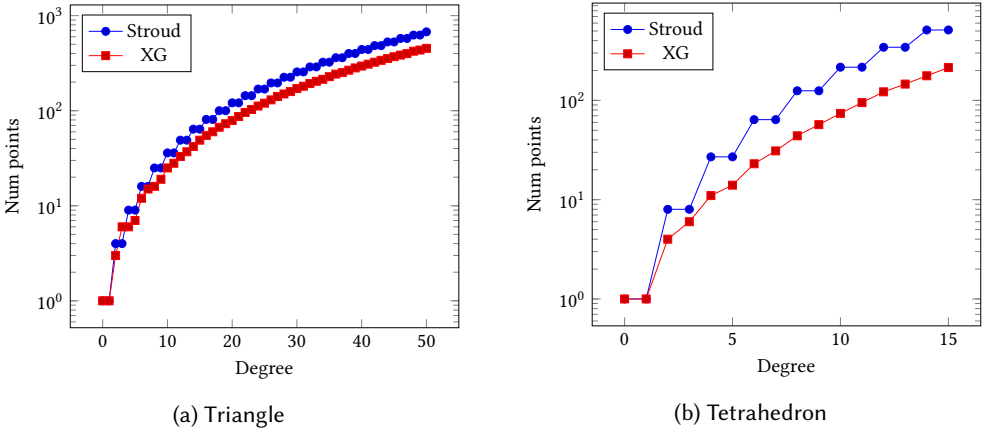

The main practical impact of these quadrature rules is the reduction of run-time in assembling high-order variational forms.
To illustrate this, we consider the elasticity operator given by the bilinear form
\begin{equation}
a(v, u) = \int_\Omega  \sigma(u) : \epsilon(v)  dx,
\end{equation}
posed over vector-valued functions.  
Here,
\begin{equation}
\epsilon(u) = \frac{1}{2} \left( \nabla u + {\left( \nabla u\right)}^\transpose\right)
\end{equation}
is the symmetric strain rate tensor, and the stress tensor is given by
\begin{equation}
\sigma(u) = \lambda (\nabla \cdot u) I + 2 \mu \epsilon(u),
\end{equation}
and $\lambda, \mu$ are the Lam\'e parameters.

\begin{figure}
\begin{subfigure}[c]{0.49\textwidth}
\centering
\begin{tikzpicture}[scale=0.75]
\begin{semilogyaxis}[xlabel={Degree}, ylabel={Time (s)}, legend pos=north west]
\addplot table [x=deg,y=action, col sep=comma] {elasticity.2d.canonical.csv};
\addlegendentry{Action (Stroud)}
\addplot table [x=deg,y=action, col sep=comma] {elasticity.2d.default.csv};
\addlegendentry{Action (XG)}
\addplot table [x=deg,y=assembly, col sep=comma] {elasticity.2d.canonical.csv};
\addlegendentry{Assembly (Stroud)}
\addplot table [x=deg,y=assembly, col sep=comma] {elasticity.2d.default.csv};
\addlegendentry{Assembly (XG)}
\end{semilogyaxis}
\end{tikzpicture}
\caption{2D\label{assembly2d}}
\end{subfigure}
\begin{subfigure}[c]{0.49\textwidth}
\centering
\begin{tikzpicture}[scale=0.75]
\begin{semilogyaxis}[xlabel={Degree}, ylabel={Time (s)}, legend pos=north west]
\addplot table [x=deg,y=action, col sep=comma] {elasticity.3d.canonical.csv};
\addlegendentry{Action (Stroud)}
\addplot table [x=deg,y=action, col sep=comma] {elasticity.3d.default.csv};
\addlegendentry{Action (XG)}
\addplot table [x=deg,y=assembly, col sep=comma] {elasticity.3d.canonical.csv};
\addlegendentry{Assembly (Stroud)}
\addplot table [x=deg,y=assembly, col sep=comma] {elasticity.3d.default.csv};
\addlegendentry{Assembly (XG)}
\end{semilogyaxis}
\end{tikzpicture}
\caption{3D\label{assembly3d}}
\end{subfigure}
\caption{Run-time to assemble the elasticity operator and its matrix-free action for various polynomial degrees with Xiao--Gimbutas and Stroud conical rules.\label{assembly}}
\end{figure}
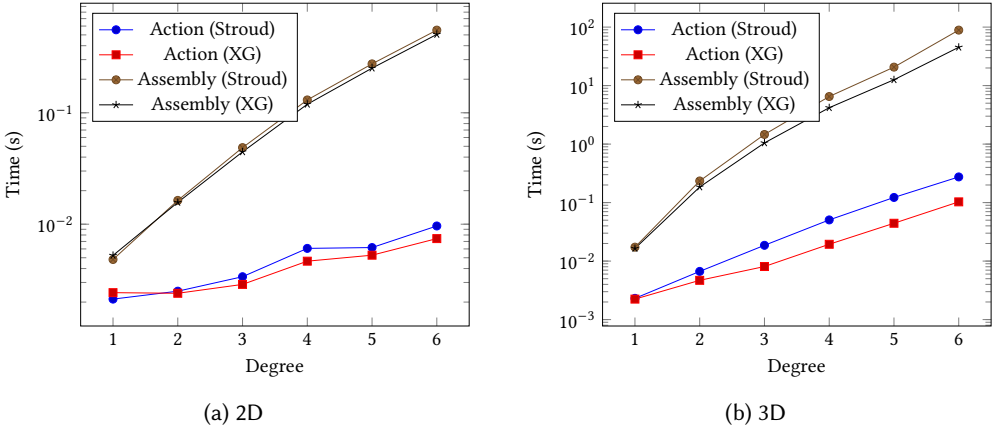

In 2D, we take $\Omega = [0, 1] \times [0, 0.2]$ divided into $64 \times 32$ rectangles, each subdivided into right triangles.
In 3D, we take $\Omega = [0, 1] \times [0, 0.2] \times [0, 0.2]$ divided into $16 \times 8 \times 8$ boxes, each subdivided into six tetrahedra.
We then measure the cost of assembling the stiffness matrix associated with $a(u, v)$ using Lagrange basis functions of degree $1 \leq k \leq 6$, integrating with either Xiao--Gimbutas or Stroud conical rules of degree $2 k$.  
We also measure the cost of computing the matrix-free action of the operator by integration:
\begin{equation}
a(\phi_i, u) = \int_\Omega \sigma(u) : \epsilon(\phi_i) dx,
\end{equation}
where $\phi_i$ is the nodal basis for the finite element space.  
These timings are shown in \cref{assembly}, where we see a modest reduction in run-time from using the Xiao--Gimbutas rule.
The savings are more pronounced in three dimensions, as the reduction in number of quadrature points relative to the Stroud rule is greater.

Effective high-order quadrature rules afford several practical advantages.
First, the cost of assembling finite element variational forms is directly proportional to the number of quadrature points, so clients immediately obtain a run-time improvement in assembling some higher-order operators.
Second, improved integration schemes reduce the cost of element instantiation and code generation for elements defined with integral moments as degrees of freedom.
Third, computing the nodal interpolant of functions with these elements requires evaluating the integral moments, which also greatly benefits from the reduction in number of quadrature points.

\section{Elements for fast diagonalization methods}
Finally, we briefly discuss the implementation of the tensor-product
elements of \citet{brubeck23} used to obtain solvers with optimal complexity in polynomial degree.
These discretize the $L^2$-de Rham complex on hexahedral cells, and 
are inspired by the fast diagonalization method (FDM) of \citet{lynch64},
to obtain sparse stiffness and mass matrices.
In one dimension, the FDM degrees of freedom for $\mathrm{CG}_p$, 
the $C^0$-continuous element of degree $p$, 
consist of point evaluation at the vertices and
integral moments against a numerically-computed polynomial
basis ${\{\hat{s}_i\}}_{i=1}^{p-1}$ for $\mathbb{P}_p([-1,1]) \cap H^1_0$ that is
orthogonal in both the $L^2$- and $H^1$-inner products,
\begin{equation}
   (\hat{s}_i, \hat{s}_j) = \delta_{ij}, \qquad (\hat{s}'_i, \hat{s}'_j) = \lambda_i \delta_{ij},
   \qquad \hat{s}_i(\pm 1) = 0, \qquad \forall\, i,j\in 1,\ldots, p-1.
\end{equation}
The interior-orthogonality of the basis functions gives
rise to sparse mass and stiffness matrices in 1D, and this sparsity also
extends to higher dimensions, rendering the high-order FDM operators
as sparse as a lowest-order discretization on a mesh with the same number
of degrees of freedom. In addition, these elements are amenable to
static condensation techniques: the interior block is diagonal,
giving rise to a sparse interface Schur complement.

The extension to $\mathrm{Q}_p$ elements is done by taking tensor products of
the one-dimensional bases. In \cref{tab:fdm-primal} we show numerical results
for solving the Poisson equation using conjugate gradients preconditioned by a
statically-condensed two-level additive Schwarz method with vertex-star patches
and a $p=1$ coarse space \citep{pavarino1993additive}.  One important observation
is that the memory requirements of this solver scale suboptimally, as the
Cholesky factorization of the patch submatrices requires $\mathcal{O}(p^{2d})$
nonzeros \citep{brubeck2022scalable}.

\begin{table}[htbp]
\centering
\caption{CG iteration counts and number of nonzero entries in the Cholesky factors 
   required to solve the primal Poisson formulation discretized with tensor-product FDM elements.\label{tab:fdm-primal} }
\csvreader[
   head to column names,
   tabular=r*2{S[table-format=9]}r*2{S[table-format=1.2e-1,scientific-notation=true,round-mode=figures,round-precision=3]},
   table head=\toprule {$p$} & {$\dim(V_{h,p})$} & {nonzeros} & {iterations} & {$L^2$ error} & {$H^1$ error} \\\midrule,
   table foot=\bottomrule
]{fdm_primal_poisson_afw.csv}{}
   {\degree{} & \dimV{} & \nonzeros{} & \iterations{} & \error{} & \errorgrad{}}
\end{table}

To overcome the elevated memory costs, one alternative is to relax the
continuity of the approximation and introduce a dual mixed formation in $\hdiv\times
L^2$. Tensor product $L^2$ elements are fully discontinuous, and tensor product
elements for $\hdiv$ only need continuity of the normal component across
faces. These are constructed 
in the usual way \citep{nedelec1980mixed, arnold2015finite, mcrae2016automated} by taking
tensor-products of the appropriate continuous $\mathrm{CG}_p$ or discontinuous
$\mathrm{DG}_{p-1}$ basis along each direction on each vector component. 
The key ingredient to construct a sparse preconditioner is to choose 
the degrees of freedom for the FDM variant for $\mathrm{DG}_{p-1}$ as integral 
moments against $\{1\} \cup {\{\hat{s}'_i\}}_{i=1}^{p-1}$, up to a scale factor.

\begin{figure}[htbp]
\begin{subfigure}[c]{0.49\textwidth}
\centering
\includegraphics[width=0.9\textwidth]{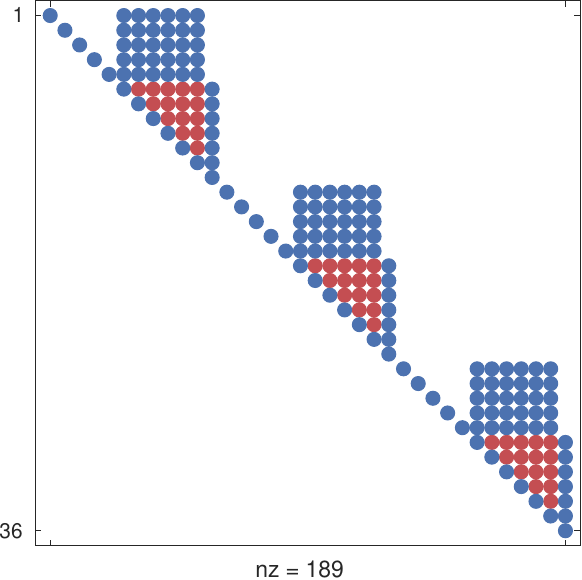}
   \caption{$\mathrm{NCF}_3(\star e)$\label{sparsity3}}
\end{subfigure}
\begin{subfigure}[c]{0.49\textwidth}
\centering
\includegraphics[width=0.9\textwidth]{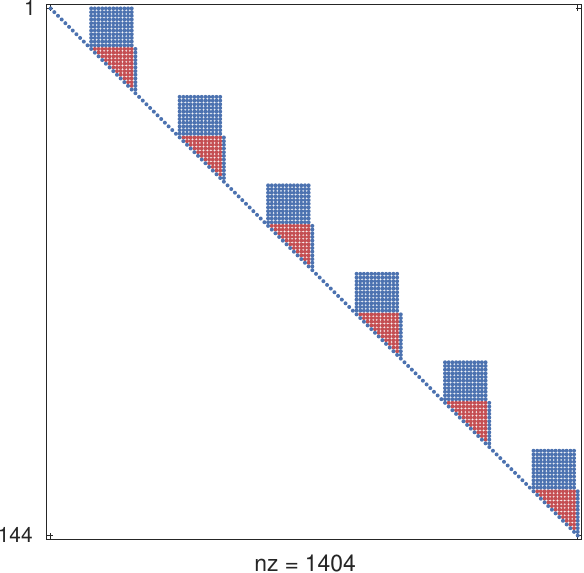}
   \caption{$\mathrm{NCF}_6(\star e)$\label{sparsity6}}
\end{subfigure}
\caption{
Sparsity pattern of the edge-star Schur complement submatrix for the $\hdiv$ Riesz-map.
Nonzeros of the upper triangular part of the Schur complement and the fill-in
induced by Cholesky factorization are shown in blue and red,
respectively.\label{fig:fdm-sparsity}}
\end{figure}

To  demonstrate how these elements are used to obtain a linear complexity
preconditioner, we consider the mixed FEM discretization of the Poisson
equation. The problem is recast as a system of first order PDEs by introducing
the flux $\sigma \coloneqq -\nabla u$ as an auxiliary unknown. The weak
formulation is to find $(\sigma, u) \in \hdiv \times L^2$ such that
\begin{alignat}{2}
  -(\sigma, \tau) + (u, \div\tau) &= {(g, \tau\cdot\mathbf{n})}_{\partial\Omega} &&\quad \forall \,\tau \in \hdiv,\\
   (\div\sigma, v)                &= (f, v) &&\quad \forall \,v \in L^2.
\end{alignat}
For our discretization, we obtain a regular hexahedral mesh $\mathcal{T}_h$ by
subdividing $\Omega$ into 4 cells along each axis. We select an inf-sup stable
pair of finite element spaces $\Sigma_{h,p} \times V_{h,p}$, with $\Sigma_{h,p}
= \mathrm{NCF}_p(\mathcal{T}_h) \subset \hdiv$ as the space of tensor-product
\Nedelec{} face elements of degree $p$ \citep{nedelec1980mixed,
arnold2015finite} with the FDM basis, and $V_{h,p} =
\mathrm{DQ}_{p-1}(\mathcal{T}_h) \subset L^2$ as the discontinuous
tensor-product polynomial space of degree $p-1$ with the standard Lagrange
basis on the Gau\ss--Legendre points.  The resulting discrete saddle-point
system is solved in an iterative matrix-free fashion with MINRES
\citep{paige75}, and preconditioned by a block-diagonal sparse auxiliary
operator arising from the augmented Lagrangian bilinear form
\begin{equation}
   a_\gamma((\tau, v), (\sigma, u)) = (\sigma, \tau) + (\div\sigma, \gamma \div\tau) + (u, \gamma^{-1} v).
\end{equation}
Here $\gamma$ is the augmented Lagrangian penalty parameter, set to a large
value $\gamma = 10^6$.  The scaled mass matrix on the $L^2$-block is
preconditioned with point-Jacobi, which results in a direct solver on Cartesian
meshes.  The $\hdiv$-block is further preconditioned with a (hybrid) two-level
domain decomposition method, extending the edge-based decomposition of
\citet{arnold2000multigrid} to high-order $\hdiv$ elements.
The coarse space is constructed from the lowest-order element ($p=1$), and is
combined multiplicatively in a V(1,1)-cycle. 
After static-condensation of the interior degrees of freedom, the edge-star
subdomains $\star e$ are constructed by gathering all face degrees of freedom
around each edge $e$ of the mesh, and combined additively in parallel.  The
submatrices corresponding to $\star e$ have a sparse Cholesky factorization in
the FDM basis with $\mathcal{O}(p^d)$ nonzeros, shown in \cref{fig:fdm-sparsity}.
In the high-order statically-condensed setting, this edge-based relaxation is more efficient than
the Hiptmair decomposition \citep{Hiptmair:1998} previously implemented in
\citet{brubeck23}, which has slightly larger edge submatrices on an $\hcurl$ auxiliary
space; this is in contrast with the lowest-order case where the same
edge-based space decompositions \citep{arnold2000multigrid} are less
efficient than the Hiptmair decomposition.

We set an absolute tolerance of $10^{-14}$ on the preconditioned norm of the
residual, and prescribe randomized source term $f$ and boundary data $g$, and
start from a zero initial guess.  The MINRES iteration counts and number of
nonzeros in the Cholesky factors are shown in \cref{tab:fdm-mixed}.  The
setup and application of the solver are optimal in both memory and
computational costs, as the number of nonzeros in the Cholesky factorization
per degree of freedom and the total number of MINRES iterations remain bounded
as $p$ increases.
\begin{table}[htbp]
\centering
\caption{MINRES iteration counts and number of nonzero entries in the Cholesky factors 
   required to solve the mixed Poisson formulation discretized with
   tensor-product FDM elements.\label{tab:fdm-mixed}}
\csvreader[
   head to column names,
   tabular=r*3{S[table-format=9]}r*2{S[table-format=1.2e-1,scientific-notation=true,round-mode=figures,round-precision=3]},
   table head=\toprule {$p$} & {$\dim(\Sigma_{h,p})$} & {$\dim(V_{h,p})$} & {nonzeros} & {iterations} & {$L^2$ error} & {$H^1$ error} \\\midrule,
   table foot=\bottomrule
]{fdm_mixed_poisson_afw.csv}{}
   {\degree{} & \dimSigma{} & \dimV{} & \nonzeros{} & \iterations{} & \error{} & \errorgrad{}}
\end{table}

\section{Conclusion}
This release paper describes several recent developments in FIAT to improve its feature set, performance, and accuracy for high-order finite element methods.
The improved implementation of recurrence relations described in \cref{sec:rec} allows improved performance and much higher degree polynomials to be used in practice.
This higher order capability motivates non-equispaced point distributions for the Lagrange basis, and \cref{sec:pts} describes a new interface to~\lstinline{recursivenodes}.
In the same vein, we have improved the implementation of degrees of freedom for $\hdiv$ and $\hcurl$ elements in \cref{sec:intdof}.
Rather than specifying normal or tangential components at points on each facet, we use the ``textbook'' degrees of freedom of integral moments of the relevant components.
This greatly reduces the error in the divergence or curl when interpolating non-polynomial expressions into the finite element spaces.
With greater emphasis on high-order discretization and integral-type degrees of freedom, providing efficient higher-order quadrature rules in FIAT is also important, and  \cref{sec:quad} describes the inclusion of Xiao--Gimbutas quadrature in FIAT.\@
Finally, we have implemented special bases for fast diagonalization methods on tensor-product cells and hope to extend these to simplicial elements in the future.

\appendix
\section{Reproducibility}\label{sec:reproducibility}
The exact version of Firedrake used, along with scripts employed for the generation of numerical data is archived on Zenodo \citep{zenodo-archive}.

\bibliographystyle{ACM-Reference-Format}
\bibliography{references}

\end{document}